\documentclass[11pt]{amsart}
\usepackage{amscd, amsfonts, amsthm, graphicx, amssymb}
\title{Groups of measure-preserving homeomorphisms of noncompact 2-manifolds} 
\author{Tatsuhiko Yagasaki}
\address{Department of Mathematics, Kyoto Institute of Technology, Matsugasaki, Sakyoku, Kyoto 606, Japan}
\email{yagasaki@ipc.kit.ac.jp}

\newtheorem{theorem}{Theorem}[section]
\newtheorem{proposition}{Proposition}[section] 
\newtheorem{corollary}{Corollary}[section] 
\newtheorem{lemma}{Lemma}[section]
\newtheorem{claim}{Claim}

\theoremstyle{definition}
\newtheorem{defi}{Definition}[section]

\def \cal {\mathcal}
\def \phi {\varphi}
\def \ds {\displaystyle}
\def \Lra {\Longrightarrow}
\def \lra {\longrightarrow}
\def \lla {\longleftarrow}
\def \e {\varepsilon}

\begin{document}
\baselineskip 6 mm

\thispagestyle{empty}

\maketitle

\begin{abstract}
Suppose $M$ is a noncompact connected 2-manifold and $\mu$ is a good Radon measure of $M$ with $\mu(\partial M) = 0$. 
Let ${\cal H}(M)$ denote the group of homeomorphisms of $M$ equipped with the compact-open topology and 
${\cal H}(M)_0$ denote the identity component of ${\cal H}(M)$. 
Let ${\cal H}(M; \mu)$ denote the subgroup of ${\cal H}(M)$ consisting of 
$\mu$-preserving homeomorphisms of $M$ and 
${\cal H}(M; \mu)_0$ denote the identity component of ${\cal H}(M; \mu)$. 
We use results of A.\,Fathi and R.\,Berlanga to show that ${\cal H}(M; \mu)_0$ is a strong deformation retract of ${\cal H}(M)_0$  and classify the topological type of ${\cal H}(M; \mu)_0$. 
%This implies that ${\cal H}(M; \mu)_0 \simeq {\Bbb S}^1$ if $M \cong {\Bbb R}^2, {\Bbb S}^1 \times {\Bbb R}^1, {\Bbb S}^1 \times [0, \infty), {\Bbb P}^2 - 1{\rm pt}$ and ${\cal H}(M; \mu)_0 \simeq \ast$ otherwise. 
\end{abstract}

%%%%%%%%%%%%%%%%%%%%%%%%%%%% Section 1 %%%%%%%%%%%%%%%%%%%%%%%%%%%%%%%%%%%%%%%%%%

\section{Introduction}

The purpose of this article is to study topological properties of 
the groups of measure-preserving homeomrophisms of noncompact 2-manifolds. 
Suppose $M$ is a connected 2-manifold and $X$ is a compact subpolyhedron of $M$ with respect to some triangulation of $M$. 
Let ${\cal H}_X(M)$ denote the group of homeomorphisms $h$ of $M$ such that $h|_X = id_X$, 
equipped with the compact-open topology,   
and let ${\cal H}_X(M)_0$ denote the connected component of $id_M$ in ${\cal H}_X(M)$. 
Suppose $\mu$ is a good Radon measure on $M$ such that $\mu({\rm Fr}\,X \cup \partial M) = 0$ (cf.\,\S 3). 
Let ${\cal H}_X(M; \mu)$ denote the subgroup of ${\cal H}_X(M)$ consisting of $\mu$-preserving homeomorphisms 
and let ${\cal H}_X(M, \mu)_0$ denote the connected component of $id_M$ in ${\cal H}_X(M, \mu)$. 

%A.\,Fathi \cite{Fa} and R.\,Berlanga \cite{Be3}
A.\,Fathi and R.\,Berlanga introduced an intermediate subgroup ${\cal H}_X(M, \mu\mbox{-end-reg})$ 
between ${\cal H}_X(M)$ and ${\cal H}_X(M, \mu)$. 
According to R.\,Berlanga \cite{Be3} $h \in {\cal H}(M)$ is said to be $\mu$-end-regular 
if $h$ preserves $\mu$-null sets and $\mu$-finite ends (see \S 3).  
Let ${\cal H}_X(M, \mu\mbox{-end-reg})$ denote 
the subgroup of ${\cal H}_X(M)$ consisting of $\mu$-end-regular homeomorphisms of $M$ 
and let ${\cal H}_X(M, \mu\mbox{-end-reg})_0$ denote the connected component of $id_M$ in ${\cal H}_X(M, \mu\mbox{-end-reg})$. 

When $M$ is compact, ${\cal H}_X(M)$ is an ANR \cite{LM}, cf.\,\cite{Ya1} and A.\,Fathi \cite{Fa} showed that 
${\cal H}(M, \mu)$ is a strong deformation retract of ${\cal H}(M, \mu\mbox{-end-reg})$ and 
the latter is homotopy dense in ${\cal H}(M)$. 
This means that ${\cal H}(M, \mu)$ is an ANR and a strong deformation retract of ${\cal H}(M)$. 
The topological characterization of $\ell_2$-manifold \cite{DT} implies that ${\cal H}(M, \mu)$ is a $\ell_2$-manifold. 

In the case where $M$ is noncompact, R. Berlanga \cite{Be1, Be2, Be3} 
extended the section theorem for the action of ${\cal H}(M)$ on the space of good Radon measures on $M$ \cite{OU, Fa} to the noncompact case, and showed that ${\cal H}(M, \mu)$ is a strong deformation retract of ${\cal H}(M, \mu\mbox{-end-reg})$. 
On the other hand, we have shown that ${\cal H}_X(M)_0$ is an ANR \cite{Ya2} and 
${\cal H}_X^{\rm PL}(M)_0$ is homotopy dense in ${\cal H}_X(M)_0$ \cite{Ya3}. 
Here ${\cal H}_X^{\rm PL}(M)_0$ is the connected component of $id_M$ in the group of PL-homeomorphisms of $M$ 
(with respect to any triangulation of $M$). 
Since we can isotope the triangulation of $M$ so that 
${\cal H}_X^{\rm PL}(M)_0 \subset {\cal H}_X(M, \mu\mbox{-end-reg})_0$ (\S 4), 
it follows that ${\cal H}_X(M, \mu\mbox{-end-reg})_0$ is also homotopy dense in ${\cal H}_X(M)_0$. 
Note that some sort of arguments on PL-triangulation is necessary to include the compact polyhedron $X$ in our statements. 
We can combine these results together to obtain the noncompact version of Fathi's results in dimension 2. 

\begin{theorem} Suppose $M$ is a connected 2-manifold and $X$ is a compact subpolyhedron of $M$ 
with respect to some triangulation of $M$. Then 
${\cal H}_X(M, \mu)_0$ is an ANR and it is a strong deformation retaract of ${\cal H}_X(M)_0$. 
\end{theorem} 

The homotopy type of ${\cal H}_X(M)_0$ has been classified in \cite{Ha, Ya2}. 
The infinite-dimensional manifold theory (cf.\,\cite{vM}) enables us 
to classify the topological type of ${\cal H}_X(M)_0$. 

\begin{corollary}
${\cal H}_X(M,\mu)_0$ is a topological $\ell_2$-manifold and ${\cal H}_X(M,\mu)_0 \cong P \times \ell_2$ when ${\cal H}_X(M)_0$ has the homotopy type of a compact polyhedron $P$. 
\end{corollary}

This paper is organized as follows.   
Section 2 is devoted to generalities on ANR's, $\ell_2$-manifolds, homeomorphism grous and ends of spaces.  
Section 3 includes fundametal facts on spaces of Radon measures. 
In Section 4 we show some properties of Radon measures which are necessary to prove Theorem 1.1.  

%In a succeeding article we study the subgroup of measure-preserving-homeomorphisms with compact support.

%%%%%%%%%%%%%%%%%%%%%%%%%%%%%%%%%%%%%%%%%%%%%%%%%%%%%%%%%%%%%%%%%%%%%%%%%%%%%%%%%%%%%%%%%%%%%%%%%%

\section{Homeomorphism groups of noncompact 2-manifolds}

\subsection{Conventions}

Throughout the paper %we use the following conventions: 
spaces are assumed to be separable and metrizable, and maps are always continuous (otherwise specified). The symbol $\cong$ indicates a homeomorphism and $\simeq$ denotes a homotopy equivalence (HE). 
The term ``strong deformation retract (or retraction)'' is abbreviated as SDR. 
When $A$ is a subset of a space $X$, the symbols ${\rm Fr}_X A$, ${\rm cl}_X A$ and ${\rm Int}_X A$ denote the frontier, closure and interior of $A$ relative to $X$. When $M$ is a manifold, $\partial = \partial M$ and ${\rm Int}\,M$ denote the boundary and interior of $M$ as a manifold. 
%When $N$ is an $n$-submanifold of an $n$-manifold $M$, we always assume that $N$ is a closed subset of $M$ and ${\rm Fr}_M N$ is an ($n-1$)-submanifold transversal to $\partial M$. Therefore we have ${\rm Int}\,N = {\rm Int}_M N \cap {\rm Int}\,M$ and ${\rm Fr}_M N \subset \partial N$. 

\subsection{ANR's and $\ell_2$-manifolds} \mbox{}

A metrizable space $X$ is called an ANR (absolute neighborhood retract) if any map $f : B \to X$ from a closed subset $B$ of a metrizable space $Y$ has an extension to a neighborhood $U$ of $B$. 
%If we can always take $U = Y$, then $X$ is called an AR (absolute retract). 
%An ANR is locally contractible and an AR is exactly a contractible ANR (cf.\,\cite{Hu}). 

%\begin{lemma} (\cite{Han}) 
%A metric space $X$ is an ANR iff for any $\varepsilon > 0$ there is an ANR $Y$ and maps $f : X \to Y$ and $g : Y \to X$ such that $gf$ is $\varepsilon$-homotopic to $id_X$. 
%\end{lemma}

\begin{defi} 
A subspace $B$ of a space $Y$ is said to be homotopy dense (HD) in $Y$ (or $B$ has the homotopy absorption property in $Y$)  
if there exists a homotopy $f_t : Y \to Y$ ($0 \leq t \leq 1$) such that $f_0 = id_Y$ and $f_t(Y) \subset B$ $(0 < t \leq 1)$. 
\end{defi}

\begin{lemma} If $B$ is HD in $Y$, then {\rm (i)} the inclusion $B \subset Y$ is a HE and {\rm (ii)} $Y$ is an ANR iff $B$ is an ANR \cite{Han}. 
\end{lemma}

%\begin{lemma} $B$ is HD in $Y$ iff 
%each point $y \in Y$ admits an open neighborhood V and
%a homotopy $\phi : V \times [0, 1] \to Y$ such that 
%$\phi_0$ is the inclusion $V \subset Y$ and $\phi_t(V) \subset B$ ($0 < t \leq 1$).  
%\end{lemma}

%Let $\ell_2$ denotes the separable Hilbert space $\{ (x_n) \in {\Bbb R}^\infty : \sum_n x_n^2 < \infty \}$. 
The symbol $\ell_2$ denotes the separable Hilbert space $\{ (x_n) \in {\Bbb R}^\infty : \sum_n x_n^2 < \infty \}$. 
An $\ell_2$-manifold is a separable metrizable space which is locally homeomorphic to $\ell_2$.
For topological groups there is a simple characterization of $\ell_2$-manifolds. 

\begin{theorem} {\rm (T.\,Dobrowolski - H.\,Toru\'nczyk \cite{DT})} 
A topological group $G$ is an $\ell_2$-manifold iff it is a separable, non locally compact, completely metrizable ANR.  
\end{theorem}

%The topological type of an $\ell_2$-manifold $M$ is determined by its homotopy type. 
%If $M$ has the homotopy type of a locally compact polyhedron $P$, then $M \cong P \times \ell_2$ (cf.\,\cite{vM}).  

\subsection{Homeomorphism groups of noncompact 2-manifolds} \mbox{} 

Suppose $Y$ is a locally connected, locally compact, separable metrizable space and $A$, $X$ are closed subsets of $Y$.   
Let ${\cal H}_X(Y, A)$ denote the group of homeomorphisms $h$ of $Y$ such that $h(A) = A$ and $h|_X = id_X$,  
equipped with the compact-open topology. ${\cal H}_X(Y, A)_0$ denotes the connected component of $id_M$ in ${\cal H}_X(Y, A)$.
It is known that ${\cal H}_X(Y, A)$ is a separable, completely metrizable, topological group.

\begin{defi} When $Y$ is a polyhedron, 
${\cal H}_X^{\rm PL}(Y, A)$ denotes the subgroup of ${\cal H}_X(Y, A)$ consisting of PL-homeomorphisms of $Y$ 
and ${\cal H}_X^{\rm PL}(Y, A)_0$ denotes the connected component of $id_M$ in ${\cal H}_X^{\rm PL}(Y, A)$. 
\end{defi} 

Every 2-manifold has a PL-triangulation. 
 %with respect to some triangulation of $M$. 

\begin{theorem} Suppose $M$ is a connected PL 2-manifold and $X$ is a compact subpolyhedron of $M$. 
Then {\rm (i)} ${\cal H}_X(M)_0$ is an ANR \cite{LM, Ya2} and {\rm (ii)} ${\cal H}_X^{\rm PL}(M)_0$ is HD in ${\cal H}_X(M)_0$ \cite{GH, Ya3}. 
\end{theorem}

\subsection{Ends of spaces}  (cf.\,\cite{Be3})

Suppose $Y$ is a connected, locally connected, locally compact, separable metrizable space. 
Let ${\cal K}(Y)$ denote the set of comapct subsets of $Y$ and for each $K \in {\cal K}(Y)$ let 
${\cal C}(Y - K)$ denote the set of connected components of $Y - K$.
%and ${\cal C}_u(Y - K)$ denote the subset of ${\cal C}(Y - K)$ consisting of 
%those which are not relatively compact in $Y$. 
%Note that $Y - K$ has at most finitely many non relatively compact connected components. 

\begin{defi} 
(i) An end of $Y$ is a function $e$ which assigns an $e(K) \in {\cal C}(Y - K)$ to each $K \in {\cal K}$ such that 
$e(K_1) \supset e(K_2)$ if $K_1 \subset K_2$. 

(ii) ${\cal E}(Y)$ denotes the set of ends of $Y$. 
The end compactification of $Y$ is the space $\overline{Y} = Y \cup {\cal E}(Y)$ 
equipped with the topology defined by the following conditions: 
\begin{itemize}
\item[(a)] $Y$ is an open subset of $\overline{Y}$
\item[(b)] the fundamental open neighborhoods of $e \in {\cal E}(Y)$ is given by 
\[ N(e, K) = e(K) \,\cup \,\{ e' \in {\cal E}(Y) \mid e'(K) = e(K)\} \hspace{5mm}  
(K \in {\cal K}(Y)). \]
\end{itemize}
\end{defi}

\noindent The space $\overline{Y}$ is compact, connected, metrizable and $Y$ is a dense open subset of $\overline{Y}$. 
(If $Y$ is compact, then ${\cal E}(Y) = \emptyset$ and $\overline{Y} = Y$.) 

For $h \in {\cal H}(Y)$ and $e \in {\cal E}(Y)$ an end \,$h(e) \in {\cal E}(Y)$ is defined by $h(e)(K) = h(e(h^{-1}(K)))$ \ $(K \in {\cal K}(Y))$. 
Each $h \in {\cal H}(Y)$ extends naturally to \ $\overline{h} \in {\cal H}(\overline{Y})$ \ by \ $\overline{h}(e) = h(e)$ \ ($e \in {\cal E}(Y)$). If $h \in {\cal H}(Y)_0$, then $h(e) = e$ ($e \in {\cal E}(Y)$). 

%%%%%%%%%%%%%%%%%%%%%%%%%%%%%%%%%%%%%%%%%%
\section{Fundamental facts on Radon measures} 

Next we recall general facts on spaces of Radon measures cf.\,\cite{Be3, Fa}. 

\subsection{Spaces of Radon measures} \mbox{} 

Suppose $Y$ is a connected, locally connected, locally compact, separable metrizable space. 
Let ${\cal B}(Y)$ denote the $\sigma$-algebra of Borel subsets of $Y$. 
A Radon measure on $Y$ is a measure $\mu$ on the measurable space $(Y, {\cal B}(Y))$ 
such that $\mu(K) < \infty$ for any compact subset $K$ of $Y$. 
Let ${\cal M}(Y)$ denote the set of Radon measures on $Y$. 

\begin{defi} The weak topology $w$ on ${\cal M}(Y)$ is the weakest topology such that the function 
\[ \mbox{$\Phi_f : {\cal M}(Y) \to {\Bbb R}$ : \ $\ds \Phi_f(\mu) = \int_Y f \,d\mu$} \]
is continuous for any continuous function \ $f : Y \to {\Bbb R}$ \ with compact support. 
The notation ${\cal M}(Y)_w$ denotes the space ${\cal M}(Y)$ equipped with the weak topology $w$. 
\end{defi}

For $\mu \in {\cal M}(Y)$ and $A \in {\cal B}(Y)$ the restriction $\mu|_A \in {\cal M}(A)$ 
is defined by $(\mu|_A)(B) = \mu(B)$ \ ($B \in {\cal B}(A)$). 

\begin{lemma} {\rm (\cite[Lemma 2.2]{Be3})}
For any closed subset $A$ of $Y$ 
the function 
\[ {\cal M}(Y)_w \to {\cal M}(A)_w : \ \mu \, \mapsto \, \mu|_A \] 
is continuous at each $\mu \in {\cal M}(Y)$ with $\mu({\rm Fr}A) = 0$. 
%\begin{itemize}
%\item[(i)] For any compact subset $K$ of $Y$ 
%the function $\Phi_K : {\cal M}(Y) \to {\Bbb R}$ : $\Phi_K(\mu) = \mu(K)$ \,is upper semicontinuous.  
%\item[(ii)] {\rm (\cite[Lemma 2.2]{Be})} For any closed subset $A$ of $Y$ 
%the function ${\cal M}(Y) \to {\cal M}(A) : \mu \mapsto \mu|_A$ is continuous at each $\mu \in {\cal M}^{{\rm Fr}\,A}(Y)$. \end{itemize}
\end{lemma}

We say that $\mu \in {\cal M}(Y)$ is good if 
$\mu(p) = 0$ for any point $p \in Y$ and $\mu(U) > 0$ for any nonempty open subset $U$ of $Y$.  
For $A \in {\cal B}(Y)$ let ${\cal M}_g^A(Y)$ denote the subset of good Radon measures $\mu$ on $Y$ with $\mu(A) = 0$.

\begin{defi} For $\mu \in {\cal M}(Y)$ \,the function $\alpha(\mu) : {\cal E}(Y) \to \{ 0, \infty \}$ is defined by 
\[ \alpha(\mu)(e) = 
\begin{cases}
\, 0 & \ \mbox{($\mu(e(K)) < \infty$ \ for some $K \in {\cal K}(Y)$)},  \\[1mm]  
\infty & \ \mbox{($\mu(e(K)) = \infty$ \ for any $K \in {\cal K}(Y)$)}.
\end{cases}
\]  
\end{defi}

We obtain the subspaces of $\mu$-finite ends and $\mu$-infinite ends, 
\[ {\cal E}_f(Y; \mu) = \{ e \in {\cal E}(Y) \mid \alpha(\mu)(e) = 0 \} \ \ \mbox{and} \ \ {\cal E}_i(Y; \mu) = \{ e \in {\cal E}(Y) \mid \alpha(\mu)(e) = \infty \}. \]  

\begin{defi} For $A, X \in {\cal B}(Y)$ and $\mu \in {\cal M}_g^A(Y)$ 
we consider the following subspaces of ${\cal M}_g^A(Y)$ : 
\begin{itemize}
\item[(i)] ${\cal M}_g^A(Y; \mu\mbox{-end-reg}) = \{ \nu \in {\cal M}_g^A(Y) \mid (a), (b), (c) \}$ : \\
\hspace{8mm} (a) $\nu(Y) = \mu(Y)$, \\
\hspace{8mm} (b) $\nu$ has the same null sets as $\mu$ (i.e., $\nu(B) = 0$ iff $\mu(B) = 0$ for any $B \in {\cal B}(Y)$), \\ 
\hspace{8mm} (c) $\alpha(\nu) = \alpha(\mu)$. 

\item[(ii)] ${\cal M}_g^A(Y, X; \mu\mbox{-end-reg}) = \{ \nu \in {\cal M}_g^A(Y; \mu) \mid (d), (e) \}$ : \\
\hspace{8mm} (d) $\nu|_X = \mu|_X$, \\
\hspace{8mm} (e) $\nu(C) = \mu(C)$ for any $C \in {\cal C}(Y - X)$. 
\end{itemize}
\end{defi}

Suppose $\mu \in {\cal M}(Y)$. 
We consider the subspace $Y \cup {\cal E}_f(Y; \mu)$ of $\overline{Y}$ and 
the space ${\cal M}(Y \cup {\cal E}_f(Y; \mu))_w$ of good Radon measures on $Y \cup {\cal E}_f(Y; \mu)$.  
Each $\nu \in {\cal M}_g(Y; \mu\mbox{-end-reg})$ has a natural extension $\overline{\nu} \in {\cal M}(Y \cup {\cal E}_f(Y; \mu))$ 
defined by $\overline{\nu}(B) = \nu(B \cap Y)$ ($B \in {\cal B}(Y \cup {\cal E}_f(Y; \mu))$). 

\begin{defi} 
The finite-end weak topology $ew$ on ${\cal M}_g(Y; \mu\mbox{-end-reg})$ is 
the weakest topology for which the following injection is continuous: 
\[ \iota : {\cal M}_g(Y; \mu\mbox{-end-reg}) \to {\cal M}(Y \cup {\cal E}_f(Y; \mu))_w \ : \ \nu \, \mapsto \, \overline{\nu}\,. \] 
The notation ${\cal M}_g(Y; \mu\mbox{-end-reg})_{ew}$ denotes 
the space ${\cal M}_g(Y; \mu\mbox{-end-reg})$ equipped with the topology $ew$. 
\end{defi}

The subspace ${\cal M}_g^A(Y, X; \mu\mbox{-end-reg})_{ew}$ has a contraction $\phi_t(\nu) = (1-t)\nu + t \mu$ ($0 \leq t \leq 1$). 

\begin{lemma} Suppose $\mu \in {\cal M}_g(Y)$, $X$ is a compact subset of $Y$ with $\mu({\rm Fr}_Y X) = 0$, $U \in {\cal C}(Y - X)$ and $A = cl_Y U$. Assume that $A$ is locally connected. Then the restriction map 
\[ r : {\cal M}_g(Y, X; \mu\mbox{-end-reg})_{ew} \ \lra \ {\cal M}_g(A; \mu|_A\mbox{-end-reg})_{ew}, \hspace{5mm} r(\nu) = \nu|_A \] 
is continuous. 
\end{lemma}

\begin{proof} We use the following notations: 
$\partial U = \{ e \in {\cal E}(Y) \mid e(X) = U\}$, $Y_1 = Y \cup {\cal E}_f(Y; \mu)$ and 
$A_1 = A \cup (\partial U \cap {\cal E}_f(Y; \mu))$. Then $A_1$ is a closed subset of $Y_1$ and the inclusion $i : A \subset Y$ induces a homeomorphism \ $\overline{i} : A \cup {\cal E}_f(A; \mu|_A) \ \cong \ A_1$. 

Consider the following commutative diagram : \\[2mm]
\hspace{20mm} 
\begin{array}[t]{l}
\begin{array}[c]{lcc}
& r & \\[-1mm]
{\cal M}_g(Y, X; \mu\mbox{-end-reg})_{ew} \hspace{6mm} & \makebox(30,2){\rightarrowfill} & 
\hspace{8mm} {\cal M}_g(A; \mu|_A\mbox{-end-reg})_{ew} \\[2mm]
\hspace{8mm}
\begin{array}[c]{c}
\mbox{\small $\cap$} \\[-2.1mm]
\hspace{1.5mm} \downarrow 
\end{array} 
\iota
& & \hspace{7mm} 
\begin{array}[c]{c}
\mbox{\small $\cap$} \\[-2.1mm]
\hspace{1.5mm} \downarrow 
\end{array} \iota_A \\[2mm] 
\end{array} \\
\begin{array}[t]{ccccc}
& & & \mbox{\small $\cong$} & \\[-3.4mm]
{\cal M}_g^{{\rm Fr}_{Y_1}A_1}(Y_1)_{w} & \ \lra \ & {\cal M}_g(A_1)_w & \ \lla \ & {\cal M}_g(A \cup {\cal E}_f(A; \mu|_A))_{w} \\[0mm]
& r_1 & & \ \overline{i}_\ast & 
\end{array}
\end{array} 
\vskip 3mm 
Here, $r(\nu) = \nu|_A$, $\iota(\nu) = \overline{\nu}$, $\iota_A(\lambda) = \overline{\lambda}$, $r_1(\lambda) = \lambda|_{A_1}$ and $\overline{i}_\ast$ is the homeomorphism induced by $\overline{i}$ (cf. \S 3.2). 
By Lemma 3.1 $r_1$ is continuous and $\iota$, $\overline{i}_\ast$ are also continuous. Thus, by Definition 3.4 the map $r$ is continuous. 
\end{proof}

\subsection{Induced measures} \mbox{}

Suppose $Y$ and $Z$ are connected, locally connected, locally compact separable metric spaces and $f: Y \to Z$ is a proper map  
($f^{-1}(K)$  is compact for any $K \in {\cal K}(Z)$). 
For $\mu \in {\cal M}(Y)$ the induced measure $f_\ast \mu \in {\cal M}(Z)$ is defined by 
$(f_\ast \mu)(C) = \mu(f^{-1}(C))$ ($C \in {\cal B}(Z)$).

\begin{lemma}
The map $f_\ast : {\cal M}(Y)_{w} \to {\cal M}(Z)_{w}$ is continous. 
\end{lemma}

Suppose $E$ is a closed subset of $Y$, $F$ is a closed subset of $Z$, $f(E) = F$ and 
$f$ maps $Y - E$ homeomorphically onto $Z - F$. 

\begin{defi} For $\nu \in {\cal M}^F(Z)$ we define $f^\ast \nu \in {\cal M}^E(Y)$ by 
\[ (f^\ast \nu)(B) = \nu(f(B - E)) \ \ (B \in {\cal B}(Y)). \] 
\end{defi}
If $\nu \in {\cal M}_g^F(Z)$ and ${\rm Int}\,E = \emptyset$, then $f^\ast \nu \in {\cal M}_g^E(Y)$.

\begin{lemma} {\rm (1)} The following maps are reciplocal homeomorphisms : 
\[ \mbox{$f_\ast : {\cal M}^E(Y)_w \to {\cal M}^F(Z)_w$ \hspace{10mm} 
$f^\ast : {\cal M}^F(Z)_w \to {\cal M}^E(Y)_w$.} \] 

{\rm (2)} If $\overline{f} : {\cal E}(Y) \to {\cal E}(Z)$ is bijective, 
then for any $\nu \in {\cal M}^F(Z)$ the following maps are reciplocal homeomorphisms : 
\[ \mbox{$f_\ast : {\cal M}^E(Y, f^\ast \nu\mbox{-end-reg})_{ew} \to {\cal M}^F(Z, \nu\mbox{-end-reg})_{ew}$ \hspace{2mm} 
$f^\ast : {\cal M}^F(Z, \nu\mbox{-end-reg})_{ew} \to {\cal M}^E(Y, f^\ast \nu\mbox{-end-reg})_{ew}$ }\]
\end{lemma}

\subsection{Groups of measure-preserving homeomorphisms} \mbox{} 

Suppose $Y$ is a connected, locally connected, locally compact, separable metrizable space and $\mu \in {\cal M}(Y)$.   

\begin{defi} Let $h \in {\cal H}(Y)$. We say that 

(i) $h$ preserves $\mu$ if $h_\ast \mu = \mu$  (i.e., $\mu(h(B)) = \mu(B)$ for any $B \in {\cal B}(Y)$),  

(ii) (\cite{Fa}) $h$ is $\mu$-biregular if $h_\ast \mu$ and $\mu$ have the same null sets 
(i.e., $\mu(h(B)) = 0$ iff $\mu(B) = 0$ for any $B \in {\cal B}(Y)$). 

(iii) (\cite{Be3}) $h$ is $\mu$-end-regular if $h$ is $\mu$-biregular and $\alpha(h_\ast \mu) = \alpha(\mu)$. 
\end{defi}

\begin{defi} Suppose $A$ and $X$ are closed subsets of $Y$. 

(i) ${\cal H}_X(Y,A; \mu)$ denotes the subgroup of ${\cal H}_X(Y,A)$ consisting of $\mu$-preserving homeomorphisms.  
${\cal H}_X(Y,A; \mu)_0$ denotes the connected components of $id_X$ in ${\cal H}_X(Y,A; \mu)$. 

(ii) ${\cal H}_X(Y,A; \mu\mbox{-end-reg})$ denotes the subgroup of ${\cal H}_X(Y,A)$ consisting of $\mu$-end-regular  homeomorphisms. 
${\cal H}_X(Y,A; \mu\mbox{-end-reg})_0$ denotes the connected components of $id_X$ in ${\cal H}_X(Y,A; \mu\mbox{-end-reg})$. 
%\item[(iii)] ${\cal H}_X(Y,A; \mu\mbox{-\,e\,-bireg})$ denotes the subgroup of ${\cal H}_X(Y,A)$ consisting of $\mu$-\,e\,-regular homeomorphisms. 
%${\cal H}_X(Y,A; \mu\mbox{-\,e\,-bireg})_0$ denotes the connected components of $id_X$ in ${\cal H}_X(Y,A; \mu\mbox{-\,e\,-reg})$. 
%\item[(iii)] 
%${\cal H}_X^\mu(Y)' = \{ h \in {\cal H}_X(Y) \mid \mu(h(C)) = \mu(C)$ \ for any connected component $C$ of $Y - X \}$. \\
%${\cal H}_X(Y; \mu\mbox{-reg})' = {\cal H}_X^\mu(Y)' \cap {\cal H}_X(Y; \mu\mbox{-reg})$. 
\end{defi} 

\begin{lemma} {\rm (1)} For $h \in {\cal H}(Y)$ 
we have $\alpha(h_\ast \mu)(\overline{h}(e)) = \alpha(\mu)(e)$ ($e \in {\cal E}(Y)$). 
In particular, if $h \in {\cal H}(Y)_0$, then $\overline{h}(e) = e$ ($e \in {\cal E}(Y)$) and $\alpha(h_\ast \mu) = \alpha(\mu)$. 

{\rm (2)} If $h \in {\cal H}_X(Y)_0$, then $h(C) = C$ for any $C \in {\cal C}(Y - X)$. 
\end{lemma}

%\begin{lemma} 
%${\cal H}_X(Y; \mu) = {\cal H}_X(Y; \mu_{Y-A})$ for any $A \in {\cal B}(X)$. 
%\end{lemma}

\subsection{Actions of homeomorphism groups on spaces of Radon measures} \mbox{} 

Suppose $Y$ is a connected, locally connected, locally compact, separable metrizable space, 
$X$ and $A$ are closed subsets of $Y$. 
The topological group ${\cal H}(Y,A)$ acts continuously on the space ${\cal M}_g^A(Y)_{w}$ by $h \cdot \nu = h_\ast \nu$. 
For each $\nu \in {\cal M}_g^A(Y)_{w}$ the subgroup ${\cal H}(Y, A; \nu)$ coincides with 
the stabilizer ${\cal H}(Y,A)_\nu$ of $\nu$ under this action.  

For $\mu \in {\cal M}_g^A(Y)$ consider the subgroup 
\[ {\cal H}_X(Y,A; \mu\mbox{-end-reg})' = 
\{ h \in {\cal H}_X(Y,A; \mu\mbox{-end-reg}) \mid \mu(h(C)) = \mu(C) \ (C \in {\cal C}(Y - X)) \}. \] 

By Lemma 3.5 (2) we have $\big({\cal H}_X(Y,A; \mu\mbox{-end-reg})'\big)_0 = {\cal H}_X(Y,A; \mu\mbox{-end-reg})_0$. 
The above action induces the continuous action of 
${\cal H}_X(Y,A; \mu\mbox{-end-reg})'$ on ${\cal M}_g^A(Y, X; \mu\mbox{-end-reg})_{ew}$. 
There exists a natural orbit map 
\[ \pi : {\cal H}_X(Y,A; \mu\mbox{-end-reg})' \to {\cal M}_g^A(Y, X; \mu\mbox{-end-reg})_{ew} \hspace{10mm} \pi(h) = h_\ast \mu. \] 
\noindent A continuous section of the orbit map $\pi$ is a map 
\[ \sigma : {\cal M}_g^A(Y, X; \mu\mbox{-end-reg})_{ew} \to {\cal H}_X(Y,A; \mu\mbox{-end-reg})_0 \]
such that $\pi \sigma = id$ \ (i.e., $\sigma(\nu)_\ast \mu = \nu$ \ ($\nu \in {\cal M}_g^A(Y, X; \mu\mbox{-end-reg})_{ew}$)). 
%If $\sigma$ is a section of $\pi$, then the map 
%\[ \sigma' : {\cal M}_g^A(Y, X; \mu)_{ew} \to {\cal H}_X(Y,A; \mu\mbox{-end-reg})_0 \hspace{5mm} 
%\sigma'(\nu) = \sigma(\nu) \sigma(\mu)^{-1} \ \ (\nu \in {\cal M}_g^A(Y, X; \mu\mbox{-end-reg})_{ew}) \] 
%is also a section of $\pi$ with $\sigma'(\mu) = id_Y$. 
%
By the notation ${\cal S}(Y, X, A, \mu)$ we mean the existence of a section of the orbit map $\pi$ for the data $(Y, X, A, \mu)$. 

\begin{lemma} Suppose ${\cal S}(Y, X, A, \mu)$ holds. Then 
\begin{itemize}
\item[(i)] $({\cal H}_{X}(Y, A; \mu\mbox{-end-reg})', {\cal H}_{X}(Y, A; \mu)) \cong {\cal H}_{X}(Y, A; \mu) \times ({\cal M}_g^A(Y, X; \mu\mbox{-end-reg})_{ew}, \{ \mu \})$,  
\item[(ii)] $({\cal H}_{X}(Y, A; \mu\mbox{-end-reg})_0, {\cal H}_{X}(Y, A; \mu)_0) \cong {\cal H}_{X}(Y, A; \mu)_0 \times ({\cal M}_g^A(Y, X; \mu\mbox{-end-reg})_{ew}, \{ \mu \})$,  
\item[(iii)] ${\cal H}_{X}(Y, A; \mu)_0$ is a SDR of ${\cal H}_{X}(Y, A; \mu\mbox{-end-reg})_0$.
\end{itemize}
\end{lemma}

\begin{proof} By the assumption the orbit map $\pi$ has a section $\sigma$. 
Replacing $\sigma(\nu)$ by $\sigma(\nu)\sigma(\mu)^{-1}$, we may assume that $\sigma(\mu) = id_Y$. 

(i) The required homeomorphism 
\[ \Phi : {\cal H}_{X}(Y, A; \mu\mbox{-end-reg})' \cong {\cal H}_{X}(Y, A; \mu) \times {\cal M}_g^A(Y, X; \mu\mbox{-end-reg})_{ew} \]
is defined by $\Phi(h) = (\sigma(h_\ast \mu)^{-1} h, h_\ast \mu)$. The inverse is given by $\Phi^{-1}(g, \nu) = \sigma(\nu) g$. 

(ii) Since $\Phi(id_Y) = (id_Y, \mu)$ and ${\cal M}_g^A(Y, X; \mu\mbox{-end-reg})_{ew}$ is connected, 
it follows that 
\[ \Phi({\cal H}_{X}(Y, A; \mu\mbox{-end-reg})_0) 
= {\cal H}_{X}(Y, A; \mu)_0 \times {\cal M}_g^A(Y, X; \mu\mbox{-end-reg})_{ew}. \] 

(iii) The singleton $\{ \mu \}$ is a SDR of ${\cal M}_g^A(Y, X; \mu\mbox{-end-reg})_{ew}$. 
\end{proof}

\begin{lemma} Suppose $Y$ and $Z$ are connected, locally connected, locally compact separable metric spaces, 
$E$ is a closed subset of $Y$ with ${\rm Int}_Y E = \emptyset$ and 
$F$ is a closed subset of $Z$ with ${\rm Int}_Z F = \emptyset$.  
Suppose $f: Y \to Z$ is a proper map, $f(E) = F$, $f$ maps $Y - E$ homeomorphically onto $Z - F$ and 
$\overline{f} : {\cal E}(Y) \to {\cal E}(Z)$ is bijective. Let $\nu \in {\cal M}^F_g(Z)$. 
We have the induced measure $f^\ast \nu \in {\cal M}^E_g(Y)$. 
Under these conditions ${\cal S}(Y, E, E, f^\ast \nu)$ implies ${\cal S}(Z, F, F, \nu)$. 
\end{lemma}

\begin{proof} Let $\mu = f^\ast \nu$. Then we have the following commutative diagram : \\[1mm]
\hspace{30mm} 
$\begin{array}[t]{ccc}
& \pi_Y & \\[-1mm]
{\cal H}_E(Y; \mu\mbox{-end-reg})_0 & \lra & {\cal M}_g^E(Y, E; \mu\mbox{-end-reg})_{ew} \\[2mm] 
\phi \Big\downarrow & & \Big\downarrow f_\ast \\[3mm] 
{\cal H}_F(Z; \nu\mbox{-end-reg})_0 & \lra & {\cal M}_g^F(Z, F; \nu\mbox{-end-reg})_{ew} \\[-1mm]
& \pi_Z & 
\end{array}$. 
\vskip 2mm 
Here, $\pi_Y$ and $\pi_Z$ are the orbit maps and $f_\ast$ is a homeomorphism with the inverse $f^\ast$ (Lemma 3.4). 
For each $h \in {\cal H}_E(Y; \mu\mbox{-end-reg})_0$ 
there exists a unique $\overline{h} \in {\cal H}_F(Z; \nu\mbox{-end-reg})_0$ with $\overline{h}f = fh$. 
The map $\phi$ is defined by $\phi(h) = \overline{h}$. 
By the assumption the orbit map $\pi_Y$ has a section $\sigma_Y$. The required section $\sigma_Z$ of the orbit map $\pi_Z$ is  defined by $\sigma_Z = \phi \sigma_Y f^\ast$. 
\end{proof}

\section{Radon measures on manifolds}
\subsection{Section theorem --- a relative version} \mbox{} 

Suppose $M$ is a connected $n$-manifold. For any $\mu \in {\cal M}_g^\partial(M)$ 
the group ${\cal H}(M; \mu\mbox{-end-reg})$ acts continuously on ${\cal M}_g^\partial(M; \mu\mbox{-end-reg})_{ew}$. 
%Then ${\cal H}(M)$ acts continuously on ${\cal M}_g^\partial(M)_{w}$ and 
%we have ${\cal H}(M)_\mu = {\cal H}(M; \mu)$ for $\mu \in {\cal M}_g^\partial(M)$. 
%For any closed subset $X$ of $M$ and $\mu \in {\cal M}_g^\partial(M)$, 
%the subgroup ${\cal H}_X(M; \mu\mbox{-end-reg})_0$ acts on ${\cal M}_g^{\partial}(M,X;\mu)_{ew}$. 

\begin{theorem} {\rm (von Neumann-Oxtoby-Ulam \cite{OU})} Suppose $M$ is a compact connected $n$-manifold. 
If $\mu, \nu \in {\cal M}_g^\partial(M)$ and $\mu(M) = \nu(M)$, then 
there exists $h \in {\cal H}_\partial(M)_0$ such that $h_\ast \mu = \nu$.
\end{theorem}

\begin{theorem} {\rm (A.\,Fathi \cite{Fa}, R.\,Berlanga \cite{Be3})} Suppose $M$ is a connected $n$-manifold. Then 
for any $\mu \in {\cal M}_g^\partial(M)$  
the orbit map $\pi : {\cal H}(M; \mu\mbox{-end-reg}) \to {\cal M}_g^\partial(M; \mu\mbox{-end-reg})_{ew}$, $\pi(h) = h_\ast \mu$ has a section 
$\sigma : {\cal M}_g^\partial(M; \mu\mbox{-end-reg})_{ew} \to {\cal H}_\partial(M; \mu\mbox{-end-reg})_0$. 
\end{theorem}

We need a relative version of this section theorem. 
%Suppose $M$ is a connected PL $n$-manifold, $\mu \in {\cal M}_g^\partial(M)$ and 
%$X$ is a compact subpolyhedron of $M$ such that $\mu({\rm Fr}\,X) = 0$.  
 
%The subgroup ${\cal H}(M; \mu\mbox{-reg})$ acts on the subspace ${\cal M}_g^\partial(M; \mu)$.  
%Suppose $M$ is an $n$-manifold, $\mu \in {\cal M}_g^\partial(M)$ and $X$ is a compact subpolyhedron of $M$ 
% such that $\mu({\rm Fr}\,X) = 0$. 

\begin{corollary} 
Suppose $M$ is a connected PL $n$-manifold, $\mu \in {\cal M}_g^\partial(M)$ and $X$ is a compact subpolyhedron  of $M$ such that $\mu({\rm Fr}\,X) = 0$. 
Then the orbit map 
\[ \pi : {\cal H}_X(M; \mu\mbox{-end-reg})' \to {\cal M}_g^\partial(M, X; \mu\mbox{-end-reg})_{ew} : \ \pi(h) = h_\ast \mu \]
has a section \ 
$\sigma : {\cal M}_g^\partial(M, X; \mu\mbox{-end-reg})_{ew} \to {\cal H}_{X \cup \partial}(M; \mu\mbox{-end-reg})_0$.  \hspace{5mm} 
%\mbox{\rm (i.e., \ \ $\lambda(\nu)_\ast \mu = \nu$ \ \ $(\nu \in {\cal M}_g^\partial(M, X; \mu)))$} \]
%such that \ \ $\lambda(\mu) = id_M$. 
\end{corollary} 

\begin{proof} 
Let $Y_i$ ($i = 1, \cdots, m$) denote the closures of connected components of $M - X$. For each $i$, 
set $\partial Y_i = ({\rm Fr}_M Y_i) \cup (Y_i \cap \partial M)$ and ${\rm Int}\,Y_i = Y_i - \partial Y_i$. 
Since %${\rm Fr}_M Y_i \subset {\rm Fr}_M X$ and 
$\partial Y_i \subset {\rm Fr}_M X \cup \partial M$, we have $\mu(\partial Y_i) = 0$, and 
by Lemma 3.2, the restriction map 
\[ \lambda_i : {\cal M}_g^\partial(M, X; \mu\mbox{-end-reg})_{ew} \to {\cal M}_g^{\partial Y_i}(Y_i; \mu|_{Y_i}\mbox{-end-reg})_{ew}, 
\ \ \lambda_i(\nu) = \nu|_{Y_i} \]
is continuous. 

Since the 2nd derived neighborhood of ${\rm Fr}_M Y_i$ in $Y_i$ is a PL-mapping cylinder neighborhood of ${\rm Fr}_M Y_i$ in $Y_i$, we can construct a connected PL $n$-manifold $N_i$ and a proper onto map $f_i : N_i \to Y_i$ such that 
$f_i(\partial N_i) = \partial Y_i$, $f_i$ maps ${\rm Int}\,N_i$ homeomorphically onto ${\rm Int}\,Y_i$ and 
$\overline{f}_i : {\cal E}(N_i) \to {\cal E}(Y_i)$ is a homeomorphism. We apply Lemma 3.7 to these data and $\mu_i = \mu|_{Y_i} \in {\cal M}_g^{\partial}(Y_i)$. By Theorem 4.2 ${\cal S}(N_i, \partial, \partial, f_i^\ast(\mu_i))$ holds and by Lemma 3.7 
${\cal S}(Y_i, \partial, \partial, \mu_i)$ also holds. Thus, we obtain a section $\sigma_i$ of the orbit map 
\[ \pi_i : {\cal H}_{\partial}(Y_i; \mu_i\mbox{-end-reg})_0 \to {\cal M}_g^{\partial}(Y_i; \mu_i\mbox{-end-reg})_{ew}, \ \ \  \pi_i(g) = g_\ast \mu_i. \]

Since $M = X \cup (\cup_i Y_i)$ and ${\rm Fr}_M\,X = \cup_i \,{\rm Fr}_M\,Y_i$,  
the required section $\sigma$ of $\pi$ is defined by 
\[ \sigma(\nu) = 
\begin{cases}
id_X & \mbox{ on } \ X \\[1mm]
\sigma_i(\lambda_i(\nu)) & \mbox{ on } \ Y_i \ \ (i = 1, \cdots, m) \hspace{10mm}
\smash{\raisebox{3.5mm}{($\nu \in {\cal M}_g^\partial(M, X; \mu\mbox{-end-reg})$).}}
\end{cases} \] 
\end{proof}

By Corollary 4.1 and Lemma 3.6 we have the following conclusion. 

\begin{corollary} Under the condition of Corolalry 4.1, for any closed subset $A$ of $\partial M$ 
\begin{itemize}
\item[(i)] $({\cal H}_{X \cup A}(M, \mu\mbox{-end-reg})_0, {\cal H}_{X \cup A}(M; \mu)_0) \cong {\cal H}_{X \cup A}(M, \mu)_0 \times ({\cal M}_g^\partial(M, X; \mu\mbox{-end-reg})_{ew}, \{ \mu \})$,  
\item[(ii)] ${\cal H}_{X \cup A}(M; \mu)_0$ is a SDR of ${\cal H}_{X \cup A}(M; \mu\mbox{-end-reg})_0$.
\end{itemize}
\end{corollary}

\subsection{PL-structures compatible with Radon measures} \mbox{} 

We show that any PL-structure can be deformed to a PL-structure compatible with a given Radon measure. 

\begin{proposition} 
Suppose $M$ is a PL $n$-manifold, $\mu \in {\cal M}_g^\partial(M)$ and $X \subset X_0$ 
are closed subpolyhedra of $M$ with $\mu({\rm Fr}_M\,X) = 0$ and $\mu(X_0 - X) = 0$. 
Then there exists a PL-structure on $M$ for which 
(i) $X$ and $X_0$ are subpolyhedra of $M$ and 
(ii) ${\cal H}^{\rm PL}_X(M)_0 \subset {\cal H}_X(M; \mu\mbox{-end-reg})$.
\end{proposition}

\begin{proof} 
The PL-structure of $M$ is given by a pair $(T, \phi)$, 
where $T$ is a simplicial complex which is a combinatorial $n$-manifold and $\phi : |T| \cong M$ is a homeomorphism.   
Since $X$ and $X_0$ are subpolyhedron of $M$, subdividing $T$ if necessary, we may assume that there are subcomplexes $S$ and $S_0$ of $T$ such that  
$X = \phi(|S|)$ and $X_0 = \phi(|S_0|)$. Let $T^{(i)}$ denote the $i$-skeleton of $T$, while $T_{(i)}$ denotes the set of $i$-simplexes of $T$.

\begin{claim} For each $i = 0, \cdots, n-1$ \\[1mm] 
\hspace*{1mm} $(\ast)_i$ there exists a PL-isotopy $f_t^{i} \in {\cal H}_{X_0 \cup \partial M}^{\rm PL}(M)_0$ such that $f_0^{i} = id_M$ and $\mu(f_1^{i}(\phi(|T^{(i)}|) - X)) = 0$. 
\end{claim}

We proceeds by the induction on $i$. 

\noindent $(\ast)_0$ : Since $\mu$ is a good measure, $\mu(\phi(|T^{(0)}|)) = 0$ and we can take $f_t^{0} = id_M$. 

\noindent $(\ast)_{i-1} \Lra (\ast)_i$ ($i = 1, \cdots, n - 1$) : \ 
Given the isotopy $f_t^{i-1}$ in $(\ast)_{i-1}$. 
Let $\psi = f_1^{i-1} \phi$ and consider the barycentric subdivision $\mbox{st}\,T$ of $T$. 
For every $\sigma \in T_{(i)}$ we put $B_\sigma = \mbox{st}\,(b(\sigma), \mbox{st}\,T)$ 
(the star of the barycenter $b(\sigma)$ of $\sigma$ in $\mbox{st}\,T$).  
Then (i) $B_\sigma$ is a PL $n$-ball, \ $\sigma \cap \partial B_\sigma = \partial \sigma$, \ 
(ii) $|T| = \cup_{\sigma \in T_{(i)}} B_\sigma$, \ $B_\sigma \cap B_\tau = \partial B_\sigma \cap \partial B_\tau$ \ ($\sigma, \tau \in T_{(i)}$, $\sigma \neq \tau$), \ 
(iii) $|T^{(i-1)}| \cap {\rm Int} B_\sigma = \emptyset$, \ $|S| \subset \cup_{\sigma \in S_{(i)}} B_\sigma$. 
The PL $n$-balls $\psi(B_\sigma)$ also have the similar properties. 

For each $\sigma \in T_{(i)}$,  
(a) if $\phi(\sigma) \not\subset X_0 \cup \partial M$, then by (i) we can find an isotopy $g_t^\sigma \in {\cal H}^{\rm PL}_{\psi(\partial B_\sigma)}(\psi(B_\sigma))$ $(t \in [0, 1])$ such that $g_0^\sigma = id_{\psi(B_\sigma)}$ and $\mu(g_1^\sigma \psi(\sigma)) = 0$, and 
(b) if $\sigma \subset X_0 \cup \partial M$, then we put $g_t^\sigma = id_{\psi(B_\sigma)}$. 
By (ii) we can define a PL-isotopy $g_t \in {\cal H}^{\rm PL}(M)_0$ by $g_t = g_t^\sigma$ on $B_\sigma$. 
Since $f_1^{i-1} = id$ on $X_0 \cup \partial M$, 
in the case (a) $\psi(B_\sigma) \cap (X_0 \cup \partial M) = f_1^{i-1}(\phi(B_\sigma) \cap (X_0 \cup \partial M)) \subset f_1^{i-1}(\phi(\partial \sigma)) \subset f_1^{i-1} \phi(\partial B_\sigma) = \psi(\partial B_\sigma)$. 
Thus we have $g_t = id$ on $X_0 \cup \partial M$. 

Define a PL-isotopy $f_t^i \in {\cal H}_{X \cup \partial M}^{\rm PL}(M)_0$ by 
$f_t^i = f_{2t}^{i-1}$ ($t \in [0, 1/2]$) and 
$f_t^i = g_{2t-1} f_1^{i-1}$ ($t \in [1/2, 1]$). 

If $\sigma \in T_{(i)}$ and $\phi(\sigma) \not\subset X_0 \cup \partial M$, 
then $\psi(\sigma) \subset \psi(B_\sigma)$ and  
$f_1^i(\phi(\sigma)) = g_1 \psi(\sigma) = g_1^\sigma\psi(\sigma)$, so $\mu(f_1^i(\phi(\sigma))) = 0$. 
Since $\mu(\partial M) = 0$ and $\mu(X_0 - X) = 0$, we have $\mu(f_1^i(\phi(|T^{(i)}|) - X)) = 0$. 
This completes the inductive step. 
%Since $|T^{(i-1)}| \cap B_\sigma \subset \partial B_\sigma$, it follows that 
%$\psi(|T^{(i-1)}|) \cap \psi(B_\sigma) \subset \psi(\partial B_\sigma)$ and $g_t = id$ on $\psi(|T^{(i-1)}|)$. 
%Thus $f_1^i(\phi(|T^{(i-1)}| - {\rm Int}\,|S|)) = g_1\psi(|T^{(i-1)}| - {\rm Int}\,|S|) = \psi(|T^{(i-1)}| - {\rm Int}\,|S|)$ has the measure 0. 
\end{proof}

\begin{claim} 
There exists an isotopy $h_t \in {\cal H}_{X_0 \cup \partial M}(M)$ such that $h_0 = id_M$ and (i), (ii) hold w.r.t. the PL-structure $\psi = h_1 \phi : |T| \cong M$. 
\end{claim}

\begin{proof} By Claim 1 there exists a PL-isotopy $f_t = f_t^{(n-1)} \in {\cal H}^{\rm PL}_{X_0 \cup \partial M}(M)_0$ \ such that \ $\mu(f_1(\phi(|T^{(n-1)}| - |S|))) = 0$.  
For any $\sigma \in T_{(n)} - S$ ($= T_{(n)} - S_0$) we have the PL $n$-ball $C_\sigma = f_1(\phi(\sigma))$. 
Since $\partial C_\sigma = f_1(\phi(\partial \sigma)) \subset f_1(\phi(|T^{(n-1)}| - {\rm Int}\,|S|)) \subset f_1(\phi(|T^{(n-1)}| - |S|)) \cup {\rm Fr}\,X$ 
and $\mu({\rm Fr}\,X) = 0$ by the assumption, 
it follows that $\mu(\partial C_\sigma) = 0$ and $\mu_\sigma := \mu|_{C_\sigma} \in {\cal M}_g^\partial(C_\sigma)$.  

Consider the Lubesgue measure $m$ on ${\Bbb R}^n$. 
The restriction of $m$ to the $n$-cube $I^n := [0, 1]^n \subset {\Bbb R}^n$ is denoted by the same symbol. 
Since any affine isomorphism of ${\Bbb R}^n$ is $m$-biregular, any PL-homeomorphism between two subpolyhedra of ${\Bbb R}^n$ is  also $m$-biregular. 

Choose a PL-homeomorphism $\alpha_\sigma : C_\sigma \cong I^n$. Then $(\alpha_\sigma)_\ast \mu_\sigma \in {\cal M}_g^\partial(I^n)$ and if we set $c_\sigma = \big((\alpha_\sigma)_\ast \mu_\sigma\big)(I^n) (= \mu(C_\sigma) > 0)$, then 
by von Neumann-Oxtoby-Ulam theorem (Theorem 4.1\,(1)) there exists an isotopy $\beta^\sigma_t \in {\cal H}_\partial(I^n)$ such that 
$\beta^\sigma_0 = id$ and $(\beta^\sigma_1)_\ast (\alpha_\sigma)_\ast \mu_\sigma = c_\sigma m$. 
We put $\gamma^\sigma_t = \beta^\sigma_t \, \alpha_\sigma : C_\sigma \cong I^n$.  

Define $g_t \in {\cal H}_{X \cup f_1\phi(|T^{(n-1)}|)}(M)_0$ by 
$g_t|_X = id_X$ and $g_t|_{C_\sigma} = (\gamma_t^\sigma)^{-1} \, \alpha_\sigma \in {\cal H}_{\partial C_\sigma}(C_\sigma)_0$ ($\sigma \in T_{(n)} - S$). 
Finally we define $h_t = g_tf_t \in {\cal H}_{X \cup \partial M}(M)_0$. 

By $M'$ and $C_\sigma'$ we denote $M$ and $C_\sigma$ which have the PL-structure given by $\psi = h_1 \phi : |T| \cong M$. Note that $C_\sigma' = \psi(\sigma)$. 
It follows that $X = \psi(|S|)$ and 
$\gamma_1^\sigma  = \beta^\sigma_1 \alpha_\sigma = (\alpha_\sigma f_1 \phi) \psi^{-1} : (C_\sigma', \mu_\sigma) \cong (I^n, c_\sigma m)$ is a measure-preserving PL-homeomorphism. 

Suppose $h \in {\cal H}^{\rm PL}_X(M')$.
By Lemma 3.5\,(1) we only have to show that $h$ is $\mu$-biregular.  
There exist subdivisions $T_1$, $T_2$ of $T$ and 
a simplicial isomorphism $k : T_1 \to T_2$ such that $h \psi = \psi |k|$.  
Take any $\tau_1 \in T_1$ with $\tau_1 \not\subset |S|$. 
There exist $\tau_2 \in T_2$ with $|k|(\tau_1) = \tau_2$ and  
$\sigma_1, \sigma_2 \in T_{(n)}$ with $\tau_1 \subset \sigma_1$, $\tau_2 \subset \sigma_2$.  
Since $|k| = id$ on $|S|$, 
it follows that $\tau_2 \not\subset |S|$ and $\sigma_1, \sigma_2 \not\subset |S|$ and we have the following diagram: 
\[ \begin{array}[t]{ccccccc}
& & & h & & & \\ 
C_{\sigma_1}' & \supset & \psi(\tau_1) & \longrightarrow & \psi(\tau_2) & \subset & C_{\sigma_2}' \\[2mm] 
\gamma_{\sigma_1} \downarrow \ \cong \ \ & & \gamma_{\sigma_1} \downarrow \ \cong \ \ & & \ \ \cong \ \downarrow \gamma_{\sigma_2} & & \ \cong \ \downarrow \gamma_{\sigma_2} \\[2mm] 
I^n & \supset & \gamma_{\sigma_1}(\psi(\tau_1)) & \longrightarrow & \gamma_{\sigma_2}(\psi(\tau_2)) & \subset & I^n
\end{array} \]
\vskip 2mm
Since $\gamma_{\sigma_2} \,h \,(\gamma_{\sigma_1})^{-1} : \gamma_{\sigma_1}(\psi(\tau_1)) \cong \gamma_{\sigma_2}(\psi(\tau_2))$ is a PL-homeomorphism, it is $m$-biregular. 
Since $\gamma_{\sigma_i} : (\psi(\tau_i), \mu) \cong (\gamma_{\sigma_i}(\psi(\tau_i)), c_\sigma m)$ ($i = 1, 2$) are measure-preserving, 
it follows that $h : \psi(\tau_1) \cong \psi(\tau_2)$ is  $\mu$-biregular. 
Since $h|_X = id_X$ is $\mu$-biregular, it follows that $h$ itself is $\mu$-biregular as required. 

This completes the proof of Proposition 4.1. 
\end{proof}

\subsection{Non locally compactness of ${\cal H}_X(M, \mu)_0$}

\begin{lemma} 
Suppose $M$ is an $n$-manifold, $X$ is a closed subset of $M$ ($X \neq M$) and $\mu \in {\cal M}_g^\partial(M)$. 
Then ${\cal H}_X(M, \mu)_0$ is not locally compact.
\end{lemma} 

\begin{proof} 
(1) First we show that ${\cal H}_\partial(J^n, m)$ is not compact, where $J^n = [-1, 1]^n$ and $m$ is the $n$-dimensional Lebesgue measure on ${\Bbb R}^n$. 
Consider the sequence of points $\mbox{\boldmath $0$}, \mbox{\boldmath $q$}_k = (0, \cdots, 0, 1/k) \in {\rm Int}\,J^n$ ($k \geq 2$). 
For each $k = 2, 3, \cdots$,  
we can take $n$-balls $D_k, E_k$ in ${\rm Int}\,J^n$ 
(of the form $[a_1, b_1] \times \cdots \times [a_n, b_n]$) 
such that $\mbox{\boldmath $0$}, \mbox{\boldmath $q$}_2 \in \partial D_k$,   
$\mbox{\boldmath $0$}, \mbox{\boldmath $q$}_k \in \partial E_k$ and $m(D_k) = m(E_k)$, and 
$h_k \in {\cal H}_\partial(J^n)$ such that $h_k(D_k) = E_k$, 
$h_k(\mbox{\boldmath $0$}) = \mbox{\boldmath $0$}$ and $h_k(\mbox{\boldmath $q$}_2) = \mbox{\boldmath $q$}_k$. 

For the homeomorphism $h_k : D_k \cong E_k$, since $(h_k)_\ast m \in {\cal M}_g^\partial(E_k)$ and $((h_k)_\ast m)(E_k) = m(D_k) = m(E_k)$, by Theorem 4.1 there exists $f_k \in {\cal H}_\partial(E_k)$ such that $(f_k)_\ast(h_k)_\ast m = m$. 
Similarly, for the homeomorphism $h_k : cl(J^n - D_k) \cong cl(J^n - E_k)$, 
since $(h_k)_\ast m \in {\cal M}_g^\partial(cl(J^n - E_k))$ and $((h_k)_\ast m)(cl(J^n - E_k)) = m(cl(J^n - D_k)) = m(cl(J^n - E_k))$, there exists $g_k \in {\cal H}_\partial(cl(J^n - E_k))$ such that $(g_k)_\ast(h_k)_\ast m = m$. 

Define $\phi_k \in {\cal H}_\partial(J^n; \mu)$ by 
$\phi_k = f_k h_k$ on $D_k$ and $\phi_k = g_k h_k$ on $cl(J^n - D_k)$. 
Since $\| \phi_k(\mbox{\boldmath $0$}) - \phi_k(\mbox{\boldmath $q$}_2) \| = 1/k \to 0$,  
any subsequence of $\phi_k$ does not converge in ${\cal H}_\partial(J^n)$ and hence   
${\cal H}_\partial(J^n; \mu)$ is not compact.  

(2) Suppose ${\cal H}_X(M, \mu)$ is locally compact. Then $id_M$ has a compact neighborhood ${\cal F}$ in ${\cal H}_X(M, \mu)$. 
There exists an $\e > 0$ such that ${\cal N}(id_M, \e) \subset {\cal F}$. 
Take an $n$-ball $B$ in $M$ such that ${\rm diam}\,B < \e$ and $\mu(\partial B) = 0$. 
Put $c = \mu(B)$. 
There is a measure-preserving homeomorphism $(J^n, c \,m) \cong (B, \mu)$ and this yields a natural closed embedding ${\cal H}_\partial(J^n, c \,m) \cong {\cal H}_\partial(B, \mu) \hookrightarrow {\cal F}$ (extending by $id$ on $M - B$). 
This contradicts the non-compactness of ${\cal H}_\partial(J^n, m)$. 
\end{proof} 

\begin{lemma} 
Suppose $M$ is an $n$-manifold, $X$ is a closed subset of $M$ and $\mu \in {\cal M}_g^\partial(M)$. 
Then ${\cal H}_X(M, \mu)$ {\rm (}or ${\cal H}_X(M, \mu)_0${\rm )} is an $\ell_2$-manifold iff it is an ANR and $X \neq M$.
\end{lemma} 

%\begin{proof}[\bf Proof of Corollary 1.1] Since $M$ is locally compact and locally connected, ${\mathcal H}(M, \mu)$ is a topological group and ${\mathcal H}_X(M, \mu)$ is a closed subgroup of ${\mathcal H}(M, \mu)$. Since $M$ is locally compact and second countable, ${\mathcal H}(M, \mu)$ is also second countable. A complete metric $\rho$ on ${\cal H}(M, \mu)$ is defined by 
%\[ \rho(f, g) = d_{\infty}(f, g) + d_{\infty}(f^{-1}, g^{-1}), \hspace{0.5cm} d_{\infty}(f, g) = \sum_{n = 1}^{\infty} \frac{1}{2^n} \sup_{x \in M_n} d(f(x), g(x)) \]
%for $f, g \in {\mathcal H}(M, \mu)$, where $d$ is a complete metric on $M$ with $d \leq 1$. Since ${\mathcal H}_X(M, \mu)_0$ is not locally compact. Finally, by Propositions 4.1, 4.2 ${\mathcal H}_X(M, \mu)_0$ is an ANR. This completes the proof.
%\end{proof}

\begin{proof}
Suppose ${\cal G} \equiv {\cal H}_X(M, \mu)$ {\rm (}or ${\cal H}_X(M, \mu)_0${\rm )} is an ANR and $X \neq M$. 
Since ${\cal H}_X(M)$ is separable and completely metrizable (cf.\,\cite{Ya2}) and 
${\cal H}_X(M, \mu)$ is a closed subgroup of ${\cal H}_X(M)$, it follows that  
${\cal G}$ is also a separable, completely metrizable topological group. 
Since $X \neq M$, by Lemma 4.1 ${\cal G}$ is not locally compact. 
Thus, by Theorem 2.1 ${\cal G}$ is an $\ell_2$-manifold. 
%we have the desired conclusion. %can conclude that ${\cal H}_X(M; \mu)_0$ is an $\ell_2$-manifold. 
\end{proof}

\section{Groups of measure-preserving homeomorphisms of 2-manifolds} 

Suppose $M$ is a connected 2-manifold and $X$ is a compact subpolyhedron of $M$ 
with respect to some triangulation of $M$. 

\begin{lemma} {\rm (1)} ${\cal H}_X(M; \mu\mbox{-end-reg})_0$ is an ANR and HD in ${\cal H}_X(M)_0$. 

{\rm (2)} ${\cal H}_X(M; \mu)_0$ is an ANR and a SDR of ${\cal H}_X(M; \mu\mbox{-end-reg})_0$.  

{\rm (3)} ${\cal H}_X(M; \mu)_0$ is a SDR of ${\cal H}_X(M)_0$. 
\end{lemma} 

\begin{proof} (1) By Propsoition 4.1 $M$ has a PL-triangulation such that 
$X$ is a subpolyhedron and ${\cal H}^{\rm PL}_X(M)_0 \subset {\cal H}_X(M; \mu\mbox{-end-reg})_0$. 
Since ${\cal H}^{\rm PL}_X(M)_0 \subset {\cal H}_X(M; \mu\mbox{-end-reg})_0 \subset {\cal H}_X(M)_0$ 
and ${\cal H}^{\rm PL}_X(M)_0$ is HD in ${\cal H}_X(M)_0$ (\cite[Theorem 3.2]{Ya3}), 
it follows that ${\cal H}_X(M; \mu\mbox{-end-reg})_0$ is also HD in ${\cal H}_X(M)_0$. 
Since ${\cal H}_X(M)_0$ is an ANR (\cite{Ya2}), by Lemma 2.1 ${\cal H}_X(M; \mu\mbox{-end-reg})_0$ is also an ANR. 

(2) By Corollary 4.2 ${\cal H}_X(M; \mu)_0$ is a SDR of ${\cal H}_X(M; \mu\mbox{-end-reg})_0$. 
By (1) ${\cal H}_X(M; \mu)_0$ is also an ANR. 

(3) Since ${\cal H}_X(M; \mu)_0$ is a closed subset of ${\cal H}_X(M)_0$, 
using the absorbing homotopy in (1) and the SDR in (2) we can easily construct a SDR of 
${\cal H}_X(M)_0$ onto ${\cal H}_X(M; \mu)_0$. 
%Since ${\cal H}_X(M)_0$ is an ANR and ${\cal H}_X(M; \mu)_0$ is an ANR closed subset of ${\cal H}_X(M)_0$ and the inclusion 
%${\cal H}_X(M; \mu)_0 \subset {\cal H}_X(M)_0$ is a homotopy equivalence. 
%Thus ${\cal H}_X(M; \mu)_0$ is a SDR of ${\cal H}_X(M)_0$. \cite[pp 30\,-\,31]{Sp}
\end{proof} 

\begin{proof}[\bf Proof of Theorem 1.1 and Corollary 1.1] \mbox{} 
The assertions follow from Lemma 5.1 and Lemma 4.2.  
\end{proof}

%%%%%%%%%%%%%%%%%%%%%%%%%%%%%%%%%%%%%%%%%% Reference %%%%%%%%%%%%%%%%%%%%%%%%%%%%%%%%%%%%%%%%%%

\end{document}